\documentstyle[amsfonts,12pt]{article}

\newtheorem{example}{Example}
\newtheorem{question}{Question}
\newtheorem{lemma}{Lemma}

\begin{document}
\title{On spaces in countable web}
\author{M. V. Matveev\\
{\small Department of Mathematics, University of California at Davis,}\\
{\small Davis, CA 95616, USA (address valid till June 30, 2000)}\\
{\small e-mail: misha$\underline{\mbox{ }}$matveev@hotmail.com}
}
\date{}
\maketitle

{\bf Abstract.}{\small 
Improving a result from the paper ``Spaces in countable web''
by Yoshikazu Yasui and Zhi-min Gao we show that Tychonoff spaces in countable discrete web
may contain closed discrete subsets of arbitrarily big cardinality.
}



\bigskip

\section{Results and discussion}

In \cite{YG}, Yoshikazu Yasui and Zhi-min Gao define a space $X$
to be {\it in countable web} provided for every open cover $\cal U$ of $X$ there 
is a countable subset $A\subset X$ such that $St(A,{\cal U})=X$.
Actually, this property was known, under several different names, long before \cite{YG}.
Thus, it was called $\omega$-star in \cite{Ikenaga}, star-Lindel\"ofness in \cite{Milena}, 
\cite{Pareek} and other papers, strong star-Lindel\"ofness in 
\cite{vDRRT} and other papers, $^*$Lindel\"ofness in \cite{MM},
\cite{SW} and other papers, countabilty of weak extent in \cite{Hod}.
A number of results, such as Theorem~3.1 are not new, too.

Further, Yasui and Gao define a space $X$
to be {\it in countable discrete web} provided for every open cover $\cal U$ of $X$ there 
is a countable, closed and discrete subset $A\subset X$ such that $St(A,{\cal U})=X$. 
This property seems to be new and interesting.
It is easily seen (and pointed out in \cite{YG}) that 
the property of being in countable discrete web is between 
being in countable web (i.e. star-Lindel\"of) and having countable extent
(i.e. all closed discrete subspaces are at most countable).
Moreower, it is sriktly between and, so to say, more close to
countable extent then to countable web.
Indeed, the examples of Tychonoff spaces in countable web which are not 
in countable discrete web are easy to find: the square of the Sorgenfrey line
is like this \cite{YG}, every pseudocompact $\Psi$-space is like this, etc.
However, the examples of spaces in countable discrete web but with uncountable extent
seem not so easy to be found.
In \cite{YG}, only a T$_1$ example (similar to example~2 in \cite{Fle})
is presented.
Here we present a ZFC Tychonoff example and a consistent normal example.
Below, $D=\{0,1\}$ is the two-point discrete space,
$\bf c$ stands for the cardinality of continuum;
the definition of the cardinal $\bf p$ can be found in \cite{vD}.

\begin{example}\label{E1}
{\em
For every cardinal $\tau$ there is a
Tychonoff space $X$ in countable discrete web with $e(X)=\tau$.
}
\end{example}

For every $\alpha<\tau$ denote by $z_\alpha$ the point in $D^\tau$
with the $\alpha$-th coordinate equal to $1$ 
and the rest of the coordinates equal to $0$.
Put $Z=\{z_\alpha:\alpha<\tau\}$ and 
$$
X=(D^\tau\times(\omega+1))\setminus((D^\tau\setminus Z)\times\{\omega\}).
$$
Then $\tilde{Z}=Z\times\{\omega\}$ is a closed discrete subspace of $X$ 
of cardinality $\tau$;
the proof that $X$ is in countable discrete web is not so straightforward;
it is presented in section~\ref{second}.

\bigskip
{\bf Remarks. 1.} 
A Tychonoff space in countable discrete web with $e(X)=\bf c$ was constructed,
also by Yan-Kui Song (\cite{Song}, Example~3.1). However, Song's construction can not
be extended to $\tau>\bf c$.
Another improvement, as compared with the Song's construction,
is that our space is separable if $\tau=\bf c$.

{\bf 2.}
The proof of $X$ in our Example~1 being in countable discrete web, presented in section 2,
is similar to the proof of Theorem~1 in \cite{Mwe}
where, for every cardinal $\tau$, a pseudocompact Tychonoff space $X$
countable web and with $e(X)\geq\tau$ is constructed.
However, the space $X$ in Example~1 is not pseudocompact.
So the following question remains open.

\begin{question}
{\em How big can be the extent of a pseudocompact Tychonoff space 
in countable discrete web?}
\end{question}

\begin{example}
($\omega_1<\bf p$) {\em
A normal space $S$ in countable discrete web with $e(S)=\bf c$.
}
\end{example}

A space $X$ is called an {\em (a)-space} \cite{Mat} provided
for every open cover $\cal U$ of $X$ and every dense subspace $Y\subset X$
there is a closed in $X$ and discrete $A\subset Y$
such that $St(A,{\cal U})=X$.
It is clear that a separable (a) space is in countable discrete web.
Now let $X\subset\Bbb R$ and let $Y_X$ denote the space 
$(X\times\{0\})\cup ({\Bbb R}\times(0,1))$ endowed with the subspace 
topology inherited from the Moore-Niemytzki plane.
By Theorem~5 from \cite{JMS}, proved by Paul Szeptycki,
$Y_X$ has property (a) as soon as $|X|<\bf p$.
So take $X\subset\Bbb R$ of cardinality $\omega_1$
and put $S=Y_X$.

Or, alternatively, take as $S$ a $\Psi$-space of cardinality $\omega_1$.
Since $|S|<\bf p$, by \cite{SV}, $S$ is an (a)-space and thus 
it is in countable discrete web.

\begin{question}{\em
Is there a ZFC example of a normal space in countable discrete web 
with uncountable extent?
}
\end{question}

\begin{question}{\em
Is there a normal space in countable web which is not in countable discrete web?
}
\end{question}

A space $X$ is called $\delta\theta$-refinable (see e.g. \cite{Burke}) if every
open cover of $X$ has an open refinement of the form $\cup\{{\cal V}_n:n\in\omega\}$
where each ${\cal V}_n$ covers $X$ and for every $x\in X$ there is $n\in\omega$
such that $|U\in{\cal V}_n:x\in U\}|<\omega$.
Aull has proved in~\cite{Aull} that every $\delta\theta$-refinable space of countable
extent is Lindel\"of.
So it is worth to note here that the spaces from Examples 1 and 2 above
are in countable discrete web, $\delta\theta$-refinable and non-Lindel\"of.

We conclude this section with one correction to \cite{YG}.
In the cited theorem from \cite{Fle} (the first lines in the section ``Preliminaries''
in \cite{YG}) starcompactness is equivalent to countable compactness 
not only in regular T$_1$ spaces but in all Hausdorff spaces as well
(in \cite{Fle} this fact was stated without proof; a proof can be found in \cite{vDRRT});
here, a space $X$ is called starcompact provided for every open cover $\cal U$ of 
$X$ there 
is a finite subset $A\subset X$ such that $St(A,{\cal U})=X$.

\section{The proof}\label{second}

Here we prove that the space $X$ form example~\ref{E1} is in countable 
discrete web.
All notation ($X$, $Z$, $\tilde{Z}$, $z_\alpha$, etc.) is like above.
We need several lemmas.
The first lemma is a weaker form of a Theorem of Fodor 
(\cite{Fodor}, see also \cite{Williams}, Theorem~3.1.5).
Let $A$ be a set and $\lambda$ a cardinal.
A {\em set mapping} of order $\lambda$ is a mapping that assigns
to each $s\in A$ a subset $f(s)\subset A$ such that $|f(s)|<\lambda$ and
$s\not\in f(s)$.
A subset $T\subset A$ is called {\em $f$-free} if $f(t)\cap T=\emptyset$
for every $t\in T$.

\begin{lemma}\label{L0}
Let $A$ be a set and $f$ a set mapping on $A$ of order $\omega$.
Then there is a countable family $\cal H$ of $f$-free subsets of $A$
such that $\cup{\cal H}=A$.
\end{lemma}

\begin{lemma}\label{L1}
For every assignment to the points $z_\alpha$ their neighbourhoods $U_\alpha$
in $D^\tau$ ($0\leq\alpha<\tau$)
there is an at most countable $S\subset D^\tau$ such that 
$S\cap U_\alpha\neq\emptyset$ for all $\alpha<\tau$ and 
$\overline{S}\cap Z=\emptyset$.
\end{lemma}

\noindent{\bf Proof:}
Without loss of generality we assume that the sets $U_\alpha$
take the form 
$U_\alpha=\{f\in D^\tau: f(\alpha)=1$ and $f(\beta)=0 \quad\forall \beta\in f(\alpha)\}$
where $f(\alpha)$ is some finite subset of ${\bf c}\setminus\{\alpha\}$.
So $f$ is a set mapping on $\tau$ of order $\omega$.
By Lemma~\ref{L0} there is a countable family $\cal H$
of $f$-free subsets of $\tau$ such that $\cup{\cal H}=\tau$.
Let $\cal H=\{H_n:n\in\omega\}$.
Without loss of generality we assume that 
$H_n\cap H_m=\emptyset$ whenever $n\neq m$ and that $|H_n|>1$ for every $n$.
Denote by $p_n$ the point in $D^\tau$ such that
$$
p_n(\alpha)=\left\{\begin{array}{l}
1\mbox{ if } \alpha\in H_n,\\
0\mbox{ otherwise.}
\end{array}\right.
$$
It is clear that $p_n\in U_\alpha$ for every $\alpha\in H_n$.
If, for some $n$ and $\alpha$, $|H_n|=\{\alpha\}$ then exactly one,
namely the $\alpha$s, coordinate of $p_n$ equals one.
In that case, redefine $p_n$ so that one more coordinate of $p_n$
equals one and still $p_n\in U_\alpha$.
After having done this we have $p_n\not\in Z$ for all $n$.

Further, put $S=\{p_n:n\in\omega\}$.
Then $S\cap U_\alpha\neq\emptyset$ for all $\alpha<\tau$.
Further, $S\cap Z=\emptyset$ and $S$ is in fact a sequence converging to the 
point with all coordinates equal to zero.
Therefore $\overline{S}\cap Z=\emptyset$.
$\Box$

\begin{lemma}\label{L2}
For every countable family $\cal U$ of nonempty open sets in $D^\tau$
there is a way to choose points $p_U\in U$ for all $U\in\cal U$ so that 
$\overline{P}\cap Z=\emptyset$ where $P=\{p_U:U\in{\cal U}\}$.
\end{lemma}

\noindent{\bf Proof:}
It is easy to see that in every nonempty open set $U\subset D^\tau$
one can pick a point, $p$, such that all but finitely many coordinates
of $p$ are equal to $1$.
Pick such a point $p_U$ in every $U\in\cal U$.
Put $A_U=\{\gamma<{\bf c}:p_U(\gamma)=0\}$ and
$A=\cup\{A_U:U\in{\cal U}\}$.
Let $\alpha<\tau$. Pick $\beta\in\tau\setminus(A\cup\{\alpha\})$.
Put $O_\alpha=\{f\in D^{\bf c}:f(\beta)=0$
and $f(\alpha)=1\}$. 
Then $O_\alpha$ is a neighbourhood of $z_\alpha$ in $D^\tau$
and $O_\alpha\cap\overline{P}=\emptyset$.
So $\overline{P}\cap Z=\emptyset$.
$\Box$

\bigskip
Now let $\cal V$ be an open cover of the space $X$ from example~\ref{E1}.
For every $n\in\omega$ put 
${\cal V}_n=\{V\cap (D^\tau\times\{n\}):V\in{\cal V}\}$.
Then ${\cal V}_n$ is an open cover of $D^\tau\times\{n\}$ 
and thus it has a finite subcover
consisting of nonempty sets, say $\widetilde{{\cal V}_n}$.

Put $\widetilde{\widetilde{{\cal V}_n}}=\{U:U\times\{n\}\in\widetilde{{\cal V}_n}\}$
and ${\cal U}=\cup\{\widetilde{\widetilde{{\cal V}_n}}:n\in\omega\}$.
Then ${\cal U}$ is a countable family of nonempty open sets in $D^\tau$.
So let $\{p_U:U\in{\cal U}\}$ be like in Lemma~\ref{L2}.
For each $n\in\omega$ put $Q_n=\{(p_U,n):U\times\{n\}\in\widetilde{{\cal V}_n}\}$
and put $Q=\cup\{Q_n:n\in\omega\}$.
It is clear that $\pi_1(Q)=P$ where 
$\pi_1:D^\tau\times(\omega+1)\to D^\tau$ is the projection of the product onto
the first factor.
Now, since $\overline{P}\cap Z=\emptyset$ we have 
$\overline{Q}\cap\tilde{Z}=\emptyset$.
Since $|Q\cap(D^\tau\times\{n\})|<\omega$
for every $n\in\omega$ it follows that $Q$ is discrete and closed in $X$.
On the other hand, $Q$ is countable and 
$St(Q,{\cal V})\supset D^\tau\times\omega$.
It remains to find another countable closed and discrete set $R\subset X$
such that $\overline{R}\cap\tilde{Z}=\emptyset$ and 
$St(R,{\cal V})\supset\tilde{Z}$.
For every $\alpha\in\tau$ choose $O_\alpha\in\cal V$ such that
$(z_\alpha,\omega)\in O_\alpha$.
Also choose $U_\alpha$ open in $D^\tau$ and $n_\alpha\in\omega$
so that $U_\alpha\times[n_\alpha,\omega]\subset O_\alpha$.
Then the sets $U_\alpha$ are like in Lemma~\ref{L1},
so let $S\subset D^\tau$ also be like in Lemma~\ref{L1}.
Enumerate $S$ as $S=\{s_k:k\in\omega\}$.
For every $n\in\omega$ put $R_n=\{(s_k,n):k\leq n\}$.
Last, put $R=\cup\{R_n:n\in\omega\}$.
It is clear that $\pi_1(R)=S$, so
$\overline{R}\cap\tilde{Z}=\emptyset$.
Again, since $|R\cap(D^\tau\times\{n\})|<\omega$
for every $n\in\omega$ it follows that $R$ is discrete
and closed in $X$. Let $\alpha\in\tau$.
Then $s_k\in U_\alpha$ for some $k\in\omega$.
Put $n=\max\{k,n_\alpha\}$.
Then $(s_k,n)\in O_\alpha$.
On the other hand, $(s_k,n)\in R$, so
$(z_\alpha,\omega)\in St(R,{\cal V})$.
$\Box$.

\bigskip
{\bf Asknowlegement.} 
The paper was written while the author was visiting the University of
California, Davis. The author expresses his gratitude to colleagues
from UC Davis for their kind hospitality.



\begin{thebibliography}{9}

\bibitem{Aull}
C. E. Aull,
{\em A generalization of a theorem of Aquaro},
Bull. Austral. Math. Soc. {\bf 9} (1973) 105-108.

\bibitem{Milena}
M. Bonanzinga,
{\em Star-Lindel\"of and absolutely star-Lindel\"of spaces},
Q$\&$A in General Topology, {\bf 16} (1998) 79-104.

\bibitem{Burke}
D. K. Burke,
{\em Covering Properties},
Handbook of Set-theoretic Topology, Edited by K.~Kunen and J.~E.~Vaughan,
Elsevier Sci. Pub. 1984, 347-422.

\bibitem{vD}
E. K. van Douwen,
{\em The integers and topology},
Handbook of Set-theoretic Topology, Edited by K.~Kunen and J.~E.~Vaughan,
Elsevier Sci. Pub. 1984, 111-167.

\bibitem{vDRRT}
E. K. van Douwen, G. M. Reed, A. W. Roscoe and I. J. Tree,
{\em Star covering properties},
Topol. and Appl. {\bf 39 } (1991) 71-103.

\bibitem{Fle}
W. M. Fleischman,
{\em A new extension of countable compactness},
Fund. Math. {\bf 67} (1970) 1-9.

\bibitem{Fodor}
G. Fodor,
{\em Proof of a conjecture of P. Erd\"os},
Acta Sci. Math. Szeged {\bf 14} (1952) 219-227.

\bibitem{Hod}
R. E. Hodel,
{\em Combinatorial set theory and cardinal function inequalities},
Proc. Amer. Math. Soc. {\bf 111} (1991) 567-575.

\bibitem{Ikenaga}
S. Ikenaga,
{\em A class which contains Lindel\"of spaces, 
separable spaces and countably compact spaces},
Memoires of Numazu College of Technology
{\bf 18} (1983) 105-108.

\bibitem{JMS}
W. Just, M. V. Matveev and P. J. Szeptycki,
{\em Some results on property (a)},
Topol. and Appl. {\bf 100 } (2000) 67-83.

\bibitem{Mat}
M. V. Matveev,
{\em Some questions on property (a)},
Q$\&$A in General Topology, {\bf 15} (1997) 103-111.

\bibitem{Mwe}
M. V. Matveev,
{\em How weak is weak extent?},
submitted.

\bibitem{MM}
Dai MuMing,
{\em A topological space cardinality inequality involving the 
$^*$Lindel\"of number},
Acta Math. Sinica, {\bf 26} (1983) 731-735.

\bibitem{Pareek}
C. M. Pareek, 
{\em On some generalizations of countably compact spaces and 
Lindel\"of spaces},
Suppl. Rend. Circ. Mat. di Palermo, Ser II {\bf 24} (1990) 169-192.

\bibitem{Song}
Yan-Kui Song,
{\em On some questions on star covering properties},
Q$\&$A in General Topology, {\bf 18} (2000) 87-92.

\bibitem{SW}
S. H. Sun and Y. M. Wang,
{\em A strengthened topological cardinal inequality},
Bull. Austral. Math. Soc. {\bf 32} (1985) 375-378.

\bibitem{SV}
P. J. Szeptycki and J. E. Vaughan,
{\em Almost disjoint families and property (a)},
Fund. Math. {\bf 158} (1998) 229-240.

\bibitem{Williams}
N.H.  Williams
{\em Combinatorial Set Theory},
North-Holland 1977.

\bibitem{YG}
Yoshikazu Yasui and Zhi-min Gao
{\it Spaces in countable web},
Houston J. of Math. {\bf 25} (1999) 327-335.

\end{thebibliography}
\end{document}